\documentclass[10pt]{article}
\usepackage[latin1]{inputenc}
\usepackage{amsmath}
\usepackage{amsfonts}
\usepackage{amssymb}
\usepackage[american]{babel}
\usepackage{graphicx}
\usepackage{amsthm}
\usepackage{pstricks}

\begin{document}

\title{On the excessive $[m]$-index of a tree}

\author{G.Mazzuoccolo \thanks{Laboratoire G-SCOP (Grenoble-INP, CNRS), Grenoble, France. Research supported by a fellowship from the European Project ``INdAM fellowships in mathematics and/or applications for experienced researchers cofunded by Marie Curie actions'' e-mail: mazzuoccolo@unimore.it} }

\maketitle

\newtheorem{proposition}{Proposition}
\newtheorem{theorem}{Theorem}
\newtheorem{lemma}{Lemma}
\newtheorem{definition}{Definition}
\newtheorem{example}{Example}
\newtheorem{corollary}{Corollary}
\newtheorem{conjecture}{Conjecture}

\begin{abstract} \noindent 
The excessive $[m]$-index of a graph $G$, denoted by $\chi'_{[m]}(G)$, is the minimum number of matchings of size $m$ needed to cover the edge-set of $G$. We set $\chi'_{[m]}(G)=\infty$ if such a cover does not exist and we  call a graph $G$ $[m]$-coverable if its excessive $[m]$-index is finite.  
Obviously $\chi'_{[1]}(G)=|E(G)|$ and it is easy to prove for a $[2]$-coverable graph $G$ that $\chi'_{[2]}(G)=\max\{\chi'(G),\lceil |E(G)|/2 \rceil \}$ holds, where $\chi'(G)$ denotes the chromatic index of $G$. The case $m=3$ is completely solved by Cariolaro and Fu in \cite{CarFu}. In this paper we prove a general formula to compute the excessive $[4]$-index of a tree and we conjecture a possible generalization for any value of $m$. Furthermore, we prove that such a formula does not work for the excessive $[4]$-index of an arbitrary graph.
\end{abstract}

\noindent \textit{ Keywords: excessive index, matchings, trees, caterpillars. }

\section{Introduction}\label{sec:intro}

Throughout this paper, a graph $G$ always means a simple connected finite graph (without loops and parallel edges).\\
We use the standard notations $V(G),E(G)$ and $\Delta(G)$ for the vertex-set, edge-set and maximum degree of a graph $G$, respectively. Further, the degree of a vertex $x\in V(G)$ will be denoted by $\delta_G(x)$. \\
A matching $M$ of $G$ is a set of independent edges of $G$ and we call $[m]$-matching a matching of $G$ of size $m$. We define a graph $G$ to be $[m]$-coverable if each edge of $G$ belongs to an $[m]$-matching (for a complete treatment of matching theory we refer to the classical book of Lovasz and Plummer \cite{LovPlu}).\\
The main aim of this paper is the study of the excessive $[m]$-index of a graph. The excessive $[m]$-index of $G$ is defined in \cite{CarFu2} as the minimum number of $[m]$-matchings needed to cover the edge-set of $G$ and denoted by $\chi'_{[m]}(G)$. (Set $\chi'_{[m]}(G)=\infty$ if $G$ is not $[m]$-coverable.)\\
The particular case $m=\frac{|V(G)|}{2}$ is largely studied, where we are asking for a cover of the edge-set with perfect matchings ($1$-regular spanning subgraphs). In this case the excessive $[m]$-index is simply called the excessive index in \cite{BonCar} and the perfect matching index in \cite{FouVan}. The author proved in \cite{Maz1} and \cite{Maz2} how two well-known conjectures of Berge, Fulkerson and Seymour (see \cite{Ful} and \cite{Sey}) can be easily stated in terms of the excessive index. Furthermore, some general properties and the exact value of the excessive index for some relevant families of graphs is studied in various recent papers (see \cite{BelYou},\cite{MazYou}, \cite{CarRiz} and \cite{Man}).\\
For what concerns small values of $m$, trivially $\chi'_{[1]}(G)=|E(G)|$ holds and it is quite easy to prove that $\chi'_{[2]}(G)=\max\{\chi'(G), \lceil |E(G)|/2 \rceil \}$ holds for each $[2]$-coverable graph $G$, where $\chi'(G)$ denotes the chromatic index of $G$. The case $m=3$ is far from being trivial and it is completely solved by Cariolaro and Fu in \cite{CarFu} (see Theorem \ref{3index}). In this paper we prove a general formula to compute the excessive $[4]$-index of a tree. This formula gives a complete answer to a problem posed in \cite{CarFu2}. Furthermore, we exhibit a graph (not a tree) for which such kind of formula does not work. Finally, we conjecture a general formula for the excessive $[m]$-index of any tree (for all values of $m$) and for the excessive $[4]$-index of a general graph.\\

\section{Notation and terminology}

In what follows, we will say that a matching $M$ of a graph $G$ can be extended to an $[m]$-matching if there exists an $[m]$-matching of $G$ which contains $M$. A matching is said to be \textit{maximal} if it cannot be extended to a larger matching.\\	
A set $\cal F$ of $[m]$-matchings of a graph $G$ such that each edge of $G$ belongs to at least one of the $[m]$-matchings is an \textit{$[m]$-cover} of $G$. An $[m]$-cover of $G$ is minimum if its size is equal to the excessive $[m]$-index of $G$.\\
The concepts of excessive $[m]$-index and $[m]$-cover are strictly related to chromatic index and edge-colorings of a graph. Indeed the chromatic index can be defined as the minimum number of matchings, without constraints on their sizes,  needed to cover the edge-set of a graph. Furthermore, an $[m]$-cover $\cal F$ of $G$ can be viewed as a \textit{multicoloring} $\cal C$ of the edge-set of $G$. More precisely, we mean that each $[m]$-matching of $\cal F$ is a color of $\cal C$ and so each edge $e \in E(G)$ receives a number of colors equal to the number of $[m]$-matchings which contain $e$. We will say that $\cal C$ is the multicoloring associated to the $[m]$-cover $\cal F$.\\
It was useful in \cite{CarFu} and \cite{CarFu2} to consider edge colorings whose color
classes have approximately the same size. The same idea turned up to be useful in our proofs as well. 
More precisely, an edge coloring is called \textit{equalized} if it has the property that the difference between the sizes of any two color classes is at most $1$.  We shall often use (without further reference) the following lemma, due to de Werra \cite{Wer} and (independently) McDiarmid \cite{Mac}.

\begin{lemma}\label{equalized}
Let $G$ be a graph. Then $G$ has an equalized edge coloring with $\chi'(G)$ colors.
\end{lemma}

\begin{figure}[h]
\centering
\includegraphics[width=5.5cm]{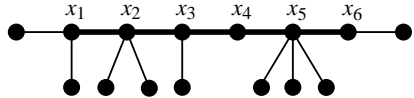} 
\caption{$CAT(1,2,1,0,3,0)$}
\label{caterpillar}
\end{figure}

Our main aim is the determination of a formula for the excessive $[4]$-index of a tree. 
A particular role in our proofs is played by caterpillars. A caterpillar is a tree in which there exists a path that contains every node of degree two or more. The (necessarily unique) path induced by the vertices of degree at least two of a caterpillar is called the \textit{spine} of the caterpillar.
We denote a caterpillar by $CAT(d_1,\ldots,d_t)$, where the path $(x_1,\ldots,x_t)$ is the spine of the caterpillar and $d_i=\delta(x_i)-2$. See Figure \ref{caterpillar} for an example.   

\section{Splitting sets}
Let $G$ be an $[m]$-coverable graph. We call a subgraph $S(m)$ of $G$ a \textit{splitting set} if 
\begin{enumerate}
\item no pair of edges of $S(m)$ can be extended to an $[m]$-matching of $G$; 
\item  $S(m)$ is not a star $K_{1,t}$ for some $t$. 
\end{enumerate}
The additional condition $2.$ is assumed in order to avoid trivial cases. When the cardinality $m$ is clear from the context, we will use the notation $S$ in place of $S(m)$. We will denote by $s(G)$ the cardinality of the largest splitting set of $G$.
The concept of a splitting set was introduced in \cite{CarFu} in order to produce a formula for the computation of the excessive $[3]$-index of an arbitrary graph.  More precisely, Cariolaro and Fu proved the following theorem.

\begin{theorem}[Cariolaro-Fu (2009)]\label{3index}
Let $G$ be a $[3]$-coverable graph. Then, $$\chi'_{[3]}(G)=\max \{ \chi'(G),\lceil |E(G)| / 3 \rceil, s(G)\}.$$
\end{theorem}

The same authors introduce in \cite{CarFu2} the notion of an $[m]$-compatible graph, that is  a graph for which $\chi'_{[m]}(G)=\max \{ \chi'(G),\lceil |E(G)| / m \rceil \}$. In particular, they prove that all trees are $[3]$-compatible.

\begin{theorem}[Cariolaro-Fu (2010)]\label{3indextrees}
Let $T$ be a $[3]$-coverable tree. Then, $$\chi'_{[3]}(T)=\max \{ \chi'(T),\lceil |E(T)|/ 3 \rceil \}.$$
\end{theorem}

\begin{figure}[h]
\centering
\includegraphics[width=11cm]{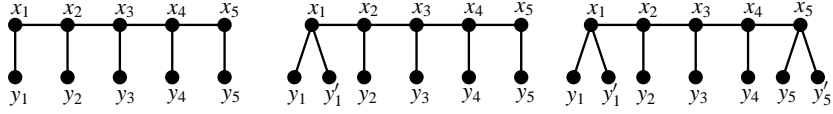} 
\caption{The non-$[4]$-compatible trees $CAT(0,1,1,1,0)$, $CAT(1,1,1,1,0)$ and $CAT(1,1,1,1,1)$.}
\label{no4comp}
\end{figure}

Furthermore, they construct a tree which is not $[5]$-compatible and they leave the determination of $[4]$-compatible trees as an open problem. In this paper, we give an answer to this question in a very strong sense by proving that there are exactly three $[4]$-coverable trees which are not $[4]$-compatible (see Figure \ref{no4comp}).

Since our main result concerns trees, we need to establish some general properties of splitting sets of an $[m]$-coverable tree.

\begin{lemma}
Let $T$ be an $[m]$-coverable tree. Let $S$ be a splitting set of $T$. Then, $\delta_{T}(x)>\delta_{S}(x), \forall x \in V(S)$. 
\end{lemma}
\textit{Proof.} Suppose $\delta_{T}(x)=\delta_S(x)$ holds for a vertex $x$ of $S$. Let $\{e_1,e_2\}$ be a $[2]$-matching of $S$ such that $e_2$ is incident to $x$ (such a $[2]$-matching does exist since $S$ is not a star). The tree $T$ is $[m]$-coverable, so $e_1$ belongs to an $[m]$-matching, say $F=\{e_1,f_2,\ldots,f_m\}$, of $T$. Note that $f_i \notin E(S)$ for $i=2,\ldots,m$, since $S$ is a splitting set. If $e_2$ is adjacent to an edge $f_i$, then $\{F \setminus \{f_i\}\} \cup \{e_2\}$ is an $[m]$-matching of $T$ containing two edges of $S$, a contradiction.
On the other hand, if $e_2$, $f_i$ are indpendent edges, then $\{F \setminus \{f_2\}\} \cup \{e_2\}$ is still an $[m]$-matching of $T$ containing two edges of $S$, a contradiction again. $\qed$ \\

\begin{lemma}
Let $T$ be an $[m]$-coverable tree. Let $S$ be a splitting set of $T$. Then, $|E(S)|\leq m$. 
\end{lemma}
\textit{Proof.} Suppose $|E(S)|\geq m+1$, then $|V(S)|\geq m+2$. Denote by $\{x_1,\ldots,x_{|V(S)|}\}$ the vertex set of $S$, and without loss of generality we suppose $[x_1,x_2]$,$[x_3,x_4]$ are two independent edges of $S$. 
Let $y_i$ be a vertex of $V(T) \setminus V(S)$ adjacent to $x_i$. The set $\{\cup_{i=5}^{|V(S)|}[x_i,y_i]\} \cup \{[x_1,x_2],[x_3,x_4]\}$ is an $[m]$-matching of $T$ containing two edges of $S$, a contradiction. $\qed$\\

\section{The excessive $[4]$-index of a tree}\label{sec:4_index_trees}

In what follows we will make large use of Lemma \ref{35covering}. This lemma is just a particular case of a more general result proved in \cite{BonMaz}. We repeat its proof in our particular context for the reader's convenience.  
 
\begin{lemma}\label{35covering}
Let ${\cal F}=\{F_1,\ldots,F_s\}$ be a set of matchings of $G$ such that 
\begin{enumerate}
 \item each edge of $G$ belongs to an element of $\cal F$,
 \item $\sum_{i=1}^s |F_i|=4s$.
\end{enumerate}
 Then, there exists a $[4]$-cover of $G$ of size $s$.
\end{lemma}
\textit{Proof.} Suppose $\cal F$ is not a $[4]$-cover of $G$, then there exist two matchings $F_i$ and $F_j$ such that $|F_i|<4<|F_j|$.
Consider the subgraph $F_i \cup F_j$, $i\neq j$, of $G$: a connected component of $F_i \cup F_j$ is either a path or a cycle of even length (possibly a single edge belonging to both $F_i$ and $F_j$). In the latter case the connected component of $F_i \cup F_j$ has the same number of edges as $F_i$ and $F_j$. Then, by $|F_i|<4<|F_j|$, there exists at least a connected component consisting of a path $P$ of odd length starting and finishing with edges of $F_j$.
The exchange of edges in $P$ increases $|F_i|$ by one and decreases $|F_j|$ by one. The iteration of this process furnishes a $[4]$-cover of $G$ of size $s$. $\qed$ \\

\begin{lemma}\label{4coverabletrees}
Let $T$ be an $[m]$-coverable tree. If $T'$ is a tree such that $V(T')=V(T)\cup\{y\}$ and $E(T')=E(T)\cup\{[x,y]\}$ where $x \in V(T)$, then $T'$ is $[m]$-coverable and $\chi'_{[m]}(T')\leq \chi'_{[m]}(T)+1$. 
\end{lemma}
\textit{Proof.} Let $e_1$ be an edge of $T$ incident $x$. Since $T$ is $[m]$-coverable $e_1$ belongs to at least an $[m]$-matching, namely $\{e_1,e_2,\ldots,e_m\}$, of $T$. We add the $[m]$-matching $\{[x,y],e_2,\ldots,e_m\}$ to an $[m]$-cover of $T$ thus obtaining an $[m]$-cover of $T'$ of size $\chi'_{[m]}(T)+1$, as required. $\qed$ \\

We prove in Theorem \ref{4compatibletrees} that each tree having a $[4]$-compatible subtree is itself compatible. Note that this property is false in general for an arbitrary graph, for instance the Petersen graph minus an edge is $[5]$-compatible but the Petersen graph is not. 

\begin{theorem}\label{4compatibletrees}
Let $T$ be a $[4]$-compatible tree. Each tree having $T$ as a subgraph is $[4]$-compatible. 
\end{theorem}
\textit{Proof.} Let $T'$ be a tree such that $V(T')=V(T)\cup\{y\}$ and $E(T')=E(T)\cup\{[x,y]\}$, where $x \in V(T)$. In order to prove that $T'$ is $[4]$-compatible we distinguish three cases, according to $\frac{|E(T)|}{4}$ being less than, equal to or greater than $\Delta(T)$, respectively.
In what follows $\cal F$ will always denote a minimum $[4]$-cover of $T$ and $\cal C$ the multicoloring of $T$ arising from $\cal F$.
Consider the case $\frac{|E(T)|}{4}<\Delta(T)$. Since $T$ is $[4]$-compatible it follows $\chi'_{[4]}(T)=\Delta(T)$. 
If $\delta_T(x)=\Delta(T)$, then $\Delta(T')=\Delta(T)+1$. By Lemma \ref{4coverabletrees}, we have $\chi'_{[4]}(T') \leq \chi'_{[4]}(T)+1 = \Delta(T)+1=\Delta(T')$, that is $\chi'_{[4]}(T')=\Delta(T')$. 
If $\delta_T(x)<\Delta(T)$, then $\Delta(T')=\Delta(T)$. Since $\frac{|E(T)|}{4}<\Delta(T)$, there exists at least an edge $f$ of $T$ receiving more than one color of $\cal C$. Denote by $c_1$ one of the colors of $f$. If $x$ is an end-vertex of $f$, then we can replace in the $[4]$-matching $c_1$ the edge $f$ with the edge $[x,y]$, producing a $[4]$-cover of $T'$ of size $\Delta(T')$. If  
all the edges incident $x$ have just one color of $\cal C$, then, by $\delta_T(x)<\Delta(T)=\chi'_{[4]}(T)$, there exists a $[4]$-matching $c$ not covering $x$. Remove the edge $f$ from the $[4]$-matching $c_1$ and add the edge $[x,y]$ to the $[4]$-matching $c$. That produces a cover of $T'$ satisfying the conditions of Lemma \ref{35covering}. Hence, there exists a $[4]$-cover of $T'$ of size $\Delta(T')$.\\
Consider the case $\frac{|E(T)|}{4}=\Delta(T)$. Since $T$ is $[4]$-compatible it follows $\chi'_{[4]}(T)=\Delta(T)=\lceil \frac{|E(T)|}{4} \rceil$. 
If $\delta_T(x)=\Delta(T)$, then $\Delta(T')=\Delta(T)+1$. By Lemma \ref{4coverabletrees}, we have $\chi'_{[4]}(T') \leq \chi'_{[4]}(T)+1 = \Delta(T)+1=\Delta(T')$, that is $\chi'_{[4]}(T')=\Delta(T')$. 
If $\delta_T(x)<\Delta(T)$, then $\frac{|E(T')|}{4} > \Delta(T')$. Hence, by Lemma 4 in \cite{CarFu}, $T'$ is $[4]$-compatible.\\
Finally, consider the case $\frac{|E(T)|}{4}>\Delta(T)$. Since $T$ is $[4]$-compatible it follows $\chi'_{[4]}(T)=\lceil \frac{|E(T)|}{4} \rceil$.
If $\delta_T(x)=\Delta(T)$, then $\Delta(T')=\Delta(T)+1$ and $|E(T')|=|E(T)|+1$. 
If $|E(T)| \equiv 0 \pmod 4$ holds, then $\lceil \frac{|E(T')|}{4} \rceil > \Delta(T')$. Hence, by Lemma 4 in \cite{CarFu}, $T'$ is $[4]$-compatible.\\
If $|E(T)| \not\equiv 0 \pmod 4$, then there exists an edge of $T$ having more than one color of $\cal C$ and we can use exactly the same argument already used in the case $\frac{|E(T)|}{4}<\Delta(T)$ and $\delta_T(x)<\Delta(T)$. This concludes the proof of the proposition. $\qed$\\

\begin{figure}[h]
\centering
\includegraphics[width=7cm]{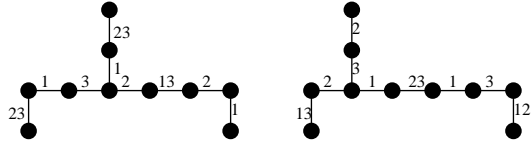} 
\caption{Two $[4]$-compatible trees of diameter $7$.}
\label{diam7a}
\end{figure}

\begin{lemma}\label{diam7}
Let $T$ be a $[4]$-coverable tree of diameter at least seven. Then $T$ is $[4]$-compatible.
\end{lemma}
\textit{Proof.} 
Let $T$ be a tree of diameter greater than $7$, that is $T$ has a path of length $8$ as subtree. Since a path of length $8$ is a $[4]$-compatible tree, it follows from Theorem \ref{4compatibletrees} that $T$ is $[4]$-compatible. From now on, we can assume the diameter of $T$ equal to $7$.\\   
Let $(x_0,x_1,\ldots,x_7)$ be a $7$-path $P$ in $T$. If $T$ is not a caterpillar then it contains one of the two trees of diameter $7$ in Figure \ref{diam7a}. Since both of them are $[4]$-compatible, it follows from Theorem \ref{4compatibletrees} that $T$ is $[4]$-compatible.\\

\begin{figure}[h]
\centering
\includegraphics[width=8cm]{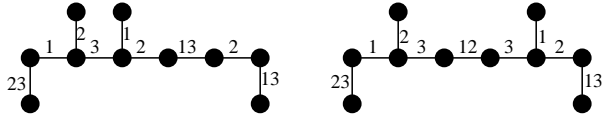} 
\caption{Two $[4]$-compatible caterpillars of diameter $7$.}
\label{diam7b}
\end{figure}

Assume that $T$ is a caterpillar $CAT(d_1,\ldots,d_6)$ with spine $(x_1,\ldots,x_6)$.
Suppose that exactly one $d_i$ is greater than $2$ and we can assume, up to symmetry, $i=1,2,3$. If $d_1>2$ (or $d_2>2$), then $[x_1,x_2]$ does not belong to a $[4]$-matching of $T$. If $d_3>2$, then $[x_3,x_4]$ does not belong to a $[4]$-matching of $T$. In any case $T$ is not $[4]$-coverable.\\
Then there are at least two $d_i$'s greater than $2$. If $d_1>2$ and $d_j>2$, for some $j$, then $T$ is not $[4]$-coverable. If $d_2>2$ and $d_4>2$, then $[x_1,x_2]$ does not belong to a $[4]$-matching. Finally, if $d_3>2$ and $d_4>2$, then $[x_3,x_4]$ does not belong to a $[4]$-matching of $T$ and, again, $T$ is not $[4]$-coverable. That leaves just two cases: $d_2>2, d_3>2$ and $d_2>2,d_5>2$. Then $T$ contains one of the two $[4]$-compatible caterpillars of Figure \ref{diam7b}, hence $T$ is $[4]$-compatible from Theorem \ref{4compatibletrees}.

In the same way one can check that if a caterpillar $T$ does not contain the two caterpillars in Figure \ref{diam7b} and it has either three or four $d_i$'s greater than $2$, then $T$ is not $[4]$-coverable. Finally, if five $d_i$'s are greater than $2$, then $T$ contains one of the caterpillars in Figure \ref{diam7b}, hence it is $[4]$-compatible from Proposition \ref{4compatibletrees}.
 $\qed$ \\

\begin{figure}[h]
\centering
\includegraphics[width=11cm]{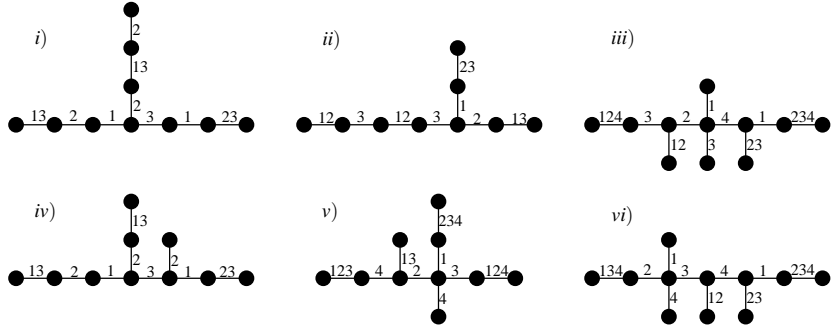} 
\caption{Six $[4]$-compatible trees.}
\label{4comp}
\end{figure}

In the following proofs we will use the fact that all trees in Figure \ref{4comp} are $[4]$-compatible (for each tree we exhibit a $[4]$-cover which uses a number of $[4]$-matchings which is equal to the maximum degree of the tree).

Furthermore, we would like to recall that the chromatic index $\chi'(T)$ of a tree is equal to its maximum degree $\Delta(T)$. 

\begin{lemma}\label{CAT01110}
Let $T$ be a tree having the caterpillar $CAT(0,1,1,1,0)$ as a subgraph. Then, $$\chi'_{[4]}(T)=\max\{\Delta(T),\lceil |E(T)|/ 4 \rceil,s(T)\}.$$
\end{lemma}
\textit{Proof.} 
First of all, we prove that the excessive $[4]$-index of the three caterpillars $CAT(0,1,1,1,0)$, $CAT(1,1,1,1,0)$ and $CAT(1,1,1,1,1)$ is equal to $4$, that is exactly the cardinality of their largest splitting set. Let $(x_1,\ldots,x_5)$ be the spine of $CAT(0,1,1,1,0)$ and let $y_i$ be the further vertex adjacent to $x_i$, for $i=1,\ldots,5$ (see Figure \ref{no4comp}). It is an easy check that  the set of edges $\{[x_1,x_2],[x_2,x_3],[x_3,x_4],[x_4,x_5]\}$ is a splitting set of size four for the caterpillar. On the other hand, a $[4]$-cover of $CAT(0,1,1,1,0)$ is given by the following four $[4]$-matchings $$\{[x_i,x_{i+1}]\} \cup \bigcup _{\begin{array}{c} j\neq i \\ j \neq i+1 \end{array}} \{ [x_j,y_j] \},$$ for $i=1,2,3,4$. The same argument can be easily adapted to obtain the same result for $CAT(1,1,1,1,0)$ and $CAT(1,1,1,1,1)$.\\
Let $x$ be an arbitrary vertex of $CAT(0,1,1,1,0)$. Now we prove that each tree $T$ obtained by adding a new vertex $z$ and the edge $[x,z]$ to $CAT(0,1,1,1,0)$ is either $[4]$-compatible or it is the caterpillar $CAT(1,1,1,1,0)$:
if $x=y_1$ ($x=y_5$), then $T$ has diameter $7$ and it is $[4]$-compatible from Lemma \ref{diam7},\\
if $x=y_2$ ($x=y_4$), then $T$ contains the graph $ii)$ in Figure \ref{4comp}, hence it is $[4]$-compatible from Theorem \ref{4compatibletrees},\\
if $x=y_3$, then $T$ contains the graph $iv)$ in Figure \ref{4comp}, hence it is $[4]$-compatible from Theorem \ref{4compatibletrees},\\
if $x=x_2$ ($x=x_4$), then $T$ contains the graph $vi)$ in Figure \ref{4comp}, hence it is $[4]$-compatible from Theorem \ref{4compatibletrees},\\
if $x=x_3$, then $T$ contains the graph $iii)$ in Figure \ref{4comp}, hence it is $[4]$-compatible from Theorem \ref{4compatibletrees}.\\
Finally, if $x=x_1$ ($x=x_5$), then $T$ is the caterpillar $CAT(1,1,1,1,0)$.
We can repeat exactly the same argument starting from $CAT(1,1,1,1,0)$ and we obtain that each tree which contains  $CAT(1,1,1,1,0)$ is either $[4]$-compatible or it is the caterpillar $CAT(1,1,1,1,1)$.
Repeating the same argument again we obtain that each tree containing $CAT(1,1,1,1,1)$ is $[4]$-compatible and the assertion follows. $\qed$\\

\begin{figure}[h]
\centering
\includegraphics[width=6cm]{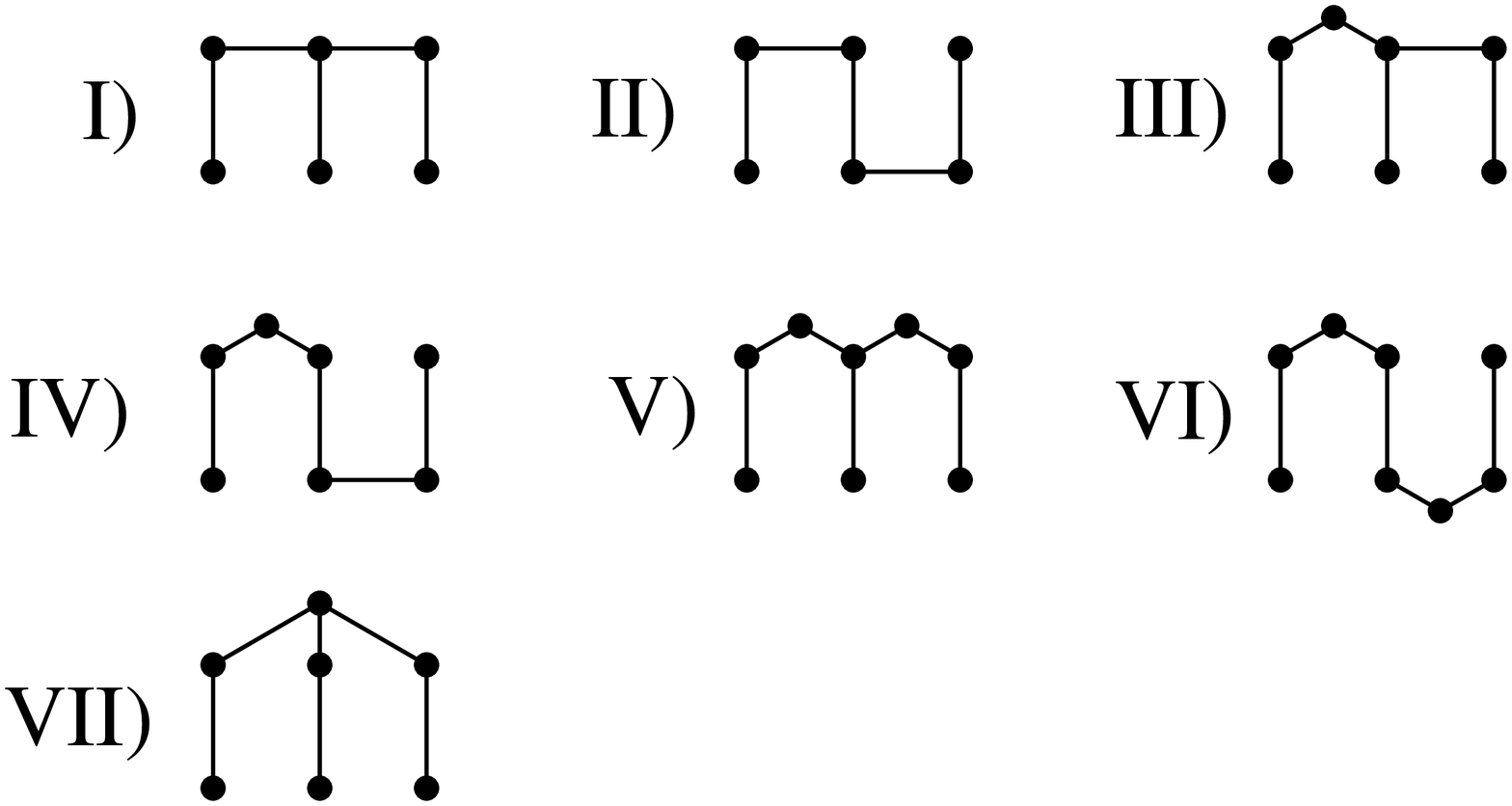} 
\caption{The seven possible cases in the proof of Theorem \ref{formula_exc_4_index_trees}.}
\label{cases}
\end{figure}

\begin{theorem}\label{formula_exc_4_index_trees}
Let $T$ be a $[4]$-coverable tree. Then, $$\chi'_{[4]}(T)=\max\{\Delta(T),\lceil |E(T)| / 4 \rceil,s(T)\}.$$
\end{theorem}
\textit{Proof.}
We can suppose $\lceil \frac{|E(T)|}{4} \rceil < \Delta(T)$, otherwise $T$ is $[4]$-compatible by Lemma 4 in \cite{CarFu} and the relation $\chi'_{[4]}(T)=\max\{\Delta(T),\lceil \frac{|E(T)|}{4} \rceil \}$ holds.\\
Suppose the existence of a maximal matching of $T$ of cardinality $2$. Hence, all edges of $T$ are incident with an edge of this $[2]$-matching and so $T$ is not $[4]$-coverable, a contradiction. On the other hand, if all maximal matchings have cardinality at least $4$, then $\chi'_{[4]}(T)=\Delta$ since each color class of any equalized $\Delta$-coloring of $T$ (see Lemma \ref{equalized}) has cardinality at most $4$ and all color classes of cardinality less than $4$ can be extended to a $[4]$-matching of $T$.
So we can assume the existence of a $[3]$-matching, say $F=\{e_1=[x_1,y_1],e_2=[x_2,y_2],e_3=[x_3,y_3]\}$, such that each edge of $T$ is incident with an edge of $F$.  
There are seven possible cases, up to symmetry, to connect the edges of $F$, as illustrated in Figure \ref{cases}. In the figure the edges of $F$ are the three vertical edges.\\

To complete our task we have to prove, for each of the seven cases presented, that either $T$ is $[4]$-compatible or $T$ contains the subgraph $CAT(0,1,1,1,0)$.

\begin{itemize}

\item \underline{Case I}: If just one of the vertices $y_i$ has degree greater than one in $T$, then the corresponding edge $[x_i,y_i]$ does not belong to a $[4]$-matching of $T$, a contradiction since $T$ is $[4]$-coverable. Then, we can assume that at least two of the vertices $y_i$ have degree greater than two in $T$, so we have three subcases:
\begin{itemize}
 \item $\delta_T(y_1)>1$, $\delta_T(y_2)>1$, $\delta_T(y_3)=1$:\\
if $\delta_T(x_1)=2$, then $[x_2,y_2]$ does not belong to a $[4]$-matching of $T$, a contradiction; if $\delta_T(x_2)=3$, then $[x_1,y_1]$ does not belong to a $[4]$-matching of $T$, a contradiction again. Hence, we have $\delta_T(x_1)>2$, $\delta_T(x_2)>3$ which implies $T$ has the tree $v)$ of Figure \ref{4comp} as a subgraph. By Theorem \ref{4compatibletrees}, $T$ is $[4]$-compatible.\\
 \item $\delta_T(y_1)>1$, $\delta_T(y_2)=1$, $\delta_T(y_3)>1$:\\
if $\delta_T(x_1)=2$ ($\delta_T(x_3)=2$), then $[x_2,x_3]$ ($[x_1,x_2]$) does not belong to a $[4]$-matching of $T$, a contradiction. Hence, $T$ contains the tree $CAT(0,1,1,1,0)$ as a subgraph, by Lemma \ref{CAT01110} the formula holds.\\
 \item $\delta_T(y_1)>1$, $\delta_T(y_2)>1$, $\delta_T(y_3)>1$:\\
if $\delta_T(x_1)=\delta_T(x_3)=2$, then $[x_2,y_2]$ does not belong to a $[4]$-matching of $T$, a contradiction. Without loss of generality we can suppose $\delta_T(x_1)>2$. Hence, $T$ contains the graph $iv)$ of Figure \ref{4comp} and $T$ is $[4]$-compatible by Theorem \ref{4compatibletrees}.
\end{itemize}

\item \underline{Case II}: If $\delta_T(y_1)=\delta_T(x_3)=1$, then $[x_2,y_2]$ does not belong to a $[4]$-matching of $T$, a contradiction. If $\delta_T(y_1)>1$ and $\delta_T(x_3)>1$, then $T$ has diameter at least $7$, so it is $[4]$-compatible by Lemma \ref{diam7}.
 Without loss of generality, we can suppose $\delta_T(y_1)>1$ and $\delta_T(x_3)=1$. 
We have $\delta_T(x_1)>2$ and $\delta_T(x_2)>2$ otherwise $[y_2,y_3]$ does not belong to a $[4]$-matching of $T$. Furthermore, $\delta_T(y_2)>2$ otherwise $[x_1,x_2]$ does not belong to a $[4]$-matching of $T$; hence, $T$ has $CAT(0,1,1,1,0)$ as subgraph and by Lemma \ref{CAT01110} the formula holds.

\item \underline{Case III}: If exactly one of the vertices $y_i$ has degree greater than one in $T$, then the corresponding edge $[x_i,y_i]$ does not belong to a $[4]$-matching of $T$, a contradiction. We have three subcases:
\begin{itemize}
\item $\delta_T(y_1)>1$, $\delta_T(y_2)>1$ and $\delta_T(y_3)=1$:\\ 
$T$ contains the graph $ii)$ of Figure \ref{4comp} and $T$ is $[4]$-compatible by Theorem \ref{4compatibletrees}.\\
\item $\delta_T(y_1)>1$, $\delta_T(y_2)=1$ and $\delta_T(y_3)>1$:\\ 
$T$ has diameter at least $7$, $T$ is $[4]$-compatible by Lemma \ref{diam7}.\\
\item $\delta_T(y_1)=1$, $\delta_T(y_2)>1$ and $\delta_T(y_3)>1$:\\ 
if $\delta_T(x_3)=2$, then the edge $[x_2,y_2]$ does not belong to a $[4]$-matching of $T$, a contradiction. Hence, $\delta_T(x_3)>2$ and $T$ contains the graph $iv)$ of Figure \ref{4comp} and  $T$ is $[4]$-compatible by Theorem \ref{4compatibletrees}.\\
\end{itemize}

\item \underline{Case IV}: If both vertices $x_3$ and $y_1$ have degree one in $T$, then the edge $[x_2,y_2]$ does not belong to a $[4]$-matching of $T$, a contradiction. Hence, $T$ has diameter at least seven and $T$ is $[4]$-compatible by Lemma \ref{diam7}.\\

\item \underline{Case V}: If just one of the vertices $y_i$ has degree greater than one in $T$, then the corresponding edge $[x_i,y_i]$ does not belong to a $[4]$-matching of $T$. Without loss of generality, we can suppose $\delta_T(y_1)>1$. Hence, $T$ has diameter at least seven and it is $[4]$-compatible by Lemma \ref{diam7}.

\item \underline{Case VI}: $T$ has diameter at least $7$ and it is $[4]$-compatible by Lemma \ref{diam7} .

\item \underline{Case VII}: If just one of the vertices $y_i$ has degree greater than one in $T$, then the corresponding edge $[x_i,y_i]$ does not belong to a $[4]$-matching of $T$, a contradiction. Without loss of generality, we can suppose $\delta_T(y_1)>1$ and $\delta_t(y_2)>1$. If $\delta_T(y_3)>1$ holds then $T$ contains the $[4]$-compatible tree $i)$ of Figure \ref{4comp} and $T$ is $[4]$-compatible by Theorem \ref{4compatibletrees}. If $\delta_T(y_3)=1$, then at least one between $x_1$ and $x_2$ has degree greater than $2$ in $T$, otherwise $T$ is not $[4]$-coverable. In both cases, $T$ contains the $[4]$-compatible tree $iv)$ of Figure \ref{4comp} and $T$ is $[4]$-compatible by Theorem \ref{4compatibletrees}. $\qed$
\end{itemize} 

The proof of the previous theorem provides also the complete list of all trees which are not $[4]$-compatible.

\begin{corollary}
The caterpillars $CAT(0,1,1,1,0)$,$CAT(1,1,1,1,0)$ and $CAT(1,1,1,1,1)$ are the unique trees which are not $[4]$-compatible. 
\end{corollary}

\section{Final remarks and conjectures}

\begin{figure}[h]
\centering
\includegraphics[width=4cm]{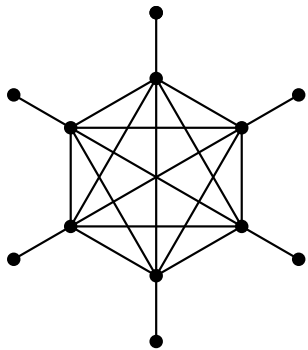} 
\caption{A graph for which the formula in Theorem \ref{formula_exc_4_index_trees} does not hold.}
\label{4index_8}
\end{figure}

Note that Theorem \ref{formula_exc_4_index_trees} cannot be generalized to the entire class of $[4]$-coverable graphs. 
Consider the graph $G$ in Figure \ref{4index_8}. It has $\Delta(G)=\chi'(G)=6$, $\lceil |E(G)|/4 \rceil = \lceil 21/4 \rceil = 6 $ and $s(G)=3$ (the largest splitting set is a circuit of length $3$). On the other hand, each $[4]$-matching of $G$ contains at most two edges in the clique of size six of $G$. Hence, $\chi'_{[4]}(G) \geq \lceil 15 / 2 \rceil = 8 $. It is possible to verify that $G$ admits a $[4]$-cover of size $8$, that is $\chi'_{[4]}(G)=8$.\\

The same idea can be used to prove that an analogous formula does not work to compute the excessive $[5]$-index of all trees. Consider the caterpillar $K=CAT(0,1,1,1,1,1,0)$. Each $[5]$-matching shares at most two edges with the spine of $K$. Since the spine has seven edges, we have $\chi'_{[5]}(K) \geq \lceil 7 / 2 \rceil = 4$, but it is easy to check that $\max \{ \chi'(K), \lceil |E(K)| / 5 \rceil, s(K) \} = \max \{3,3,1\} = 3$.\\

These examples suggest a natural generalization of the concept of a splitting set. Let $G$ be an $[m]$-coverable graph. We call a subgraph $S^t(m)$ of $G$ a $t$-splitting set if $i)$ no $[t+1]$-matching of $S^t(m)$ can be extended to an $[m]$-matching of $G$; $ii)$ $S^t(m)$ has a $[t+1]$-matching. We will denote by $s^t(G)$ the cardinality of the largest $t$-splitting set of $G$ (for a fixed $m$).\\

It is straightforward that a $1$-splitting set is exactly a splitting set previously defined and $s(G)=s^1(G)$ holds.
With this new concept in our hands, we propose the following conjectures:

\begin{conjecture}
Let $T$ be an $[m]$-coverable tree. Then, $$\chi'_{[m]}(T)=\max\{\chi'(T),\lceil |E(T)| / m \rceil,s^1(T),\lceil s^2(T) / 2 \rceil,\ldots,\lceil s^{m-1}(T) / (m-1) \rceil\}.$$ 
\end{conjecture}

\begin{conjecture}
Let $G$ be a $[4]$-coverable graph. Then, $$\chi'_{[4]}(G)=\max\{\chi'(G),\lceil |E(G)| / 4 \rceil,s^1(G),\lceil s^2(G) / 2 \rceil, \lceil s^{3}(G) / 3 \rceil\}.$$ 
\end{conjecture}

We would like to stress that the natural generalization of the second conjecture to the case $m=5$ is not possible due to the Petersen graph.\\
 
A possible way to attack the first conjecture is the idea that $t$-splitting sets of a tree should be not very large (with respect to $m$).
Theorem \ref{4compatibletrees} can be easily generalized to arbitrary values of $m$ and so the following corollary holds. 

\begin{corollary}
Let $T$ be an $[m]$-coverable tree. If $diam(T)\geq 2m$, then $T$ is $[m]$-compatible. 
\end{corollary}
Proof. The path of length $2m$ is $[m]$-compatible.

\section{Acknowledgments}
The author is grateful to one of the anonymous reviewers for his/her very precise comments and advices.

\end{document}